\newif\ifams
\amstrue

\ifams
\documentclass[runningheads,envcountsame]{amsart}
\usepackage{amsaddr}
\else
\documentclass[runningheads,envcountsame]{llncs}
\fi
\usepackage{graphicx}
\usepackage{hyperref}

\usepackage{wrapfig}

\RequirePackage{amsmath,amssymb}
\ifams
\RequirePackage{amsthm}
\fi

\usepackage{tikz}
\usetikzlibrary{cd,calc}
\RequirePackage{hyperref}
\hypersetup{
  colorlinks=true,       
  linkcolor=blue, 
  urlcolor=blue,
  citecolor=blue, 
}
\RequirePackage{xspace}
\RequirePackage{color}

\ifams
\newtheorem{proposition}{Proposition}
\newtheorem{lemma}{Lemma}
\newtheorem{theorem}{Theorem}

\theoremstyle{definition}
\newtheorem{definition}{Definition}

\theoremstyle{remark}
\newtheorem{example}{Example}
\newtheorem{remark}{Remark}
\newenvironment{myproof}[1][]{
  \ifx&#1&
  \begin{proof}
    \else
    \begin{proof}[Proof of #1]
      \fi}{\end{proof}}
\else
\newenvironment{myproof}[1][]{
  \ifx&#1&
  \begin{proof}
  \else
  \begin{proof}[#1]
  \fi
  }{\hfill\qed\end{proof}}
\fi

\newcommand{\I}{\mathbb{I}}
\newcommand{\In}{\I_{n}}
\newcommand{\Qj}{Q_{\vee}}

\newcommand{\QjI}[1][]{Q_{\vee}(\mathbb{I}_{#1})}
\newcommand{\QjIn}{\QjI[n]}
\newcommand{\QjL}{Q_{\vee}(L)}

\RequirePackage{txfonts}

\newcommand{\set}[1]{\{\,#1\,\}}
\newcommand{\rto}[1][\;\;]{\xrightarrow{#1}}

\newcommand{\cont}{continuous\xspace}
\newcommand{\jc}{join-\cont}
\newcommand{\Jc}{Join-\cont}
\newcommand{\mc}{meet-\cont}
\newcommand{\END}{\end{document}}
\newcommand{\card}{\mathrm{card}}

\newcommand{\eqdef}{:=}
\newcommand{\Img}{Image}
\newcommand{\Image}{\Img}

\newcommand{\IMG}{J}
\newcommand{\INCR}{M}

\newcommand{\pthNoNETurns}[2]{
  \draw[black,thin] (0,0) grid (#1,#1);
  \foreach \x/\y[remember=\x as \precx,remember=\y as \precy] in {#2}
  {
    \ifdefined\precx\else\def\precx{1}\def\precy{0}\fi
    \pgfmathsetmacro\xp{\x-1};
    \coordinate (point) at (\xp,\y);
    \node at ([shift={(0.5,-0.5)}] point) {};
    \path[very  thick,black] 
    (\precx -1,\precy) edge (\xp,\precy)
    (\xp,\precy) edge (\xp,\y) ;
    \pgfmathsetmacro\precx{\xp};
    \pgfmathsetmacro\precy{\y};
  }
  \pgfmathsetmacro\previousx{\precx - 1};
  \path[very thick,black] 
  (\previousx,\precy) edge (#1,\precy)
  (#1,\precy) edge (#1,#1) ;
}

\newcommand{\diagonal}[1]{
  \path[very  thick,dotted,black]
  (0,0) edge (#1,#1);
}
\newcommand{\upDiagonal}[1]{
  \path[very  thick,dotted,black]
  (0,0.4) edge (#1-0.4,#1);
}

\newcounter{mynote}

\newcommand{\len}[2][]{|#2|_{#1}}
\newcommand{\lx}[1]{\len[x]{#1}}
\newcommand{\ly}[1]{\len[y]{#1}}

\newcommand{\blockX}[2]{\mathbf{bl}^{#1,x}_{#2}}
\newcommand{\blockY}[2]{\mathbf{bl}^{#1,y}_{#2}}
\newcommand{\blockXno}[2]{\mathbf{blno}^{#1,x}_{#2}}
\newcommand{\BLOCKXNO}[1]{\mathbf{blno}^{#1,x}}

\newcommand{\blockYno}[2]{\mathbf{blno}^{#1,y}_{#2}}

\newcommand{\mone}{\underline{1}}

\newcommand{\emmentaler}{emmentaler\xspace}

\newcommand{\cam}{closed under arbitrary meets\xspace}
\newcommand{\caj}{closed under arbitrary joins\xspace}
\newcommand{\cajam}{closed under arbitrary joins and meets\xspace}

\renewcommand{\o}{\mathbf{int}}
\renewcommand{\j}{\mathbf{cl}}

\newcommand{\E}{\mathcal{E}}
\newcommand{\Ef}{\E_{f}}
\newcommand{\jE}{\j_{\E}}
\newcommand{\jEf}{\j_{\Ef}}
\newcommand{\oE}{\o_{\E}}
\newcommand{\oEf}{\o_{\Ef}}
\newcommand{\fE}{f_{\E}}
\newcommand{\fEf}{f_{\Ef}}
\newcommand{\gE}{g_{\E}}
\newcommand{\EfE}{\E_{\fE}}
\newcommand{\radj}{right adjoint\xspace}
\newcommand{\ladj}{left adjoint\xspace}

\renewcommand{\O}{\mathcal{O}}
\newcommand{\NE}{North-East\xspace}
\newcommand{\NETurn}{North-East turn\xspace}
\newcommand{\ENTurn}{East-North turn\xspace}

\newcommand{\NStep}{North step\xspace}
\newcommand{\ESteps}{East steps\xspace}
\newcommand{\EStep}{East step\xspace}

\newcommand{\uzigzag}{upper zigzag\xspace}
\newcommand{\OEIS}[1]{\href{https://oeis.org/#1}{#1}}

\newcommand{\map}{map\xspace}
\newcommand{\fib}{f}
\newcommand{\LU}{Laradji and Umar\xspace}

\begin{document}
\ifams
\title{On discrete idempotent paths}
\author{Luigi
  Santocanale}
\address{
  Laboratoire d'Informatique et des Syst\`emes, \\
  UMR 7020, Aix-Marseille
  Universit\'e, CNRS }
\email{luigi.santocanale@lis-lab.fr}
\thanks{Partially supported
  by the ``LIA LYSM AMU CNRS ECM INdAM'' and by the ``LIA LIRCO''}

\else
\input{llncsFrontMatter}
\fi
\maketitle              
\begin{abstract}
  The set of discrete lattice paths from $(0,0)$ to $(n,n)$ with North
  and East steps (i.e. words $w \in \set{x,y}^{\ast}$ such that
  $\len[x]{w} = \len[y]{w} = n$) has a canonical monoid structure
  inherited from the bijection with the set of \jc \map{s} from the
  chain $\set{0,1,\ldots ,n}$ to itself.
  We explicitly describe this monoid structure and, relying on a
  general characterization of idempotent \jc \map{s} from a complete
  lattice to itself, we characterize idempotent paths
  as \uzigzag paths.  We argue that these paths are counted by the odd
  Fibonacci numbers. Our method yields a geometric/combinatorial proof
  of counting results, due to Howie and to \LU, for idempotents in
  monoids of monotone endo\map{s} on finite chains.

  \ifams
  \else
  \fi
\end{abstract}

\section{Introduction}

Discrete lattice paths from $(0,0)$ to $(n,m)$ with
North and East steps have a standard representation as words
$w \in \set{x,y}^{\ast}$ such that $\len[x]{w} = n$ and
$\len[y]{w} = m$. The set $P(n,m)$ of these paths, with the dominance
ordering, is a distributive lattice (and therefore of a Heyting
algebra), see e.g. \cite{BB94,FePi05,Ferrari16,Muhle17}. A simple
proof that the dominance ordering is a lattice relies on the bijective
correspondence between these paths and monotone \map{s} from the
chain $\set{1,\ldots ,n}$ to the chain $\set{0,1,\ldots ,m}$, see e.g.
\cite{birkhoff,BB94}.  In turn, these \map{s} bijectively correspond
to \jc \map{s} from $\set{0,1,\ldots ,n}$ to $\set{0,1,\ldots ,m}$
(those order preserving \map{s} that sends $0$ to $0$).  \Jc
\map{s} from a complete lattice to itself form, when given the
pointwise ordering, a complete lattice in which composition
distributes with joins. This kind of algebraic structure combining a
monoid operation with a lattice structure is called a quantale
\cite{rosenthal1990} or (roughly speaking) a residuated lattice
\cite{GJKO}.  Therefore, the aforementioned bijection also witnesses a
richer structure for $P(n,n)$, that of a quantale
and of a residuated lattice. 
The set $P(n,n)$ is actually a \emph{star-autonomous} quantale or, as
a residuated lattice, \emph{involutive}, see
\cite{luigisPreprint}. 

A main aim of this paper is to draw attention to the interplay between
the algebraic and enumerative combinatorics of paths and these
algebraic structures (lattices, Heyting algebras, quantales,
residuated lattices) that, curiously, are all related to logic.  We
focus in this paper on the monoid structure that corresponds under the
bijection to function composition---which, from a logical perspective,
can be understood as a sort of non-commutative conjunction. In the
literature, the monoid structure appears to be less known than the
lattice structure. 
A notable exception is the work \cite{LaradjiUmar2016} where a
different kind of lattice paths, related to Delannoy paths, are
considered so to represent monoids of injective order-preserving
partial transformations on chains.

We explicitly describe the monoid structure of $P(n,n)$ and
characterize those paths that are idempotents.  Our characterization
relies on a general characterization of idempotent \jc \map{s} from a
complete lattice to itself.  When the complete lattice is the chain
$\set{0,1,\ldots ,n}$, this characterization yields a description of
idempotent paths as those paths whose all \NETurn{s} are above the
line $y = x + \frac{1}{2}$ and whose all \ENTurn{s} are below this
line.  We call these paths \emph{\uzigzag}.
We use this characterization to provide a geometric/combinatorial
proof that \uzigzag paths in $P(n,n)$ are counted by the odd Fibonacci
numbers $\fib_{2n+1}$.  Simple algebraic connections among the monoid
structure on $P(n,n)$, the monoid $\O_{n}$ of order preserving \map{s}
from $\set{1,\ldots ,n}$ to itself, and the submonoid $\O_{n}^{n}$ of
$\O_{n}$ of \map{s} fixing $n$, yield a geometric/combinatorial proof
of counting results due to Howie \cite{howie71} (the number of
idempotents in $\O_{n}$ is the even Fibonacci numbers $\fib_{2n}$) and
Laradji and Umar \cite{LaradjiUmar2006} (the number of idempotents in
$\O_{n}^{n}$ is the odd Fibonacci numbers $\fib_{2n-1}$).

\section{A product on paths}
\newcommand{\aPath}{yxxxyxyyxy}
\ifams
\begin{figure}
  \else
  \begin{wrapfigure}[13]{r}[0pt]{4.2cm}
    \fi
  \centering
  \ifams
  \vskip -10pt
  \else
  \vskip -18pt
  \fi
  \includegraphics[]{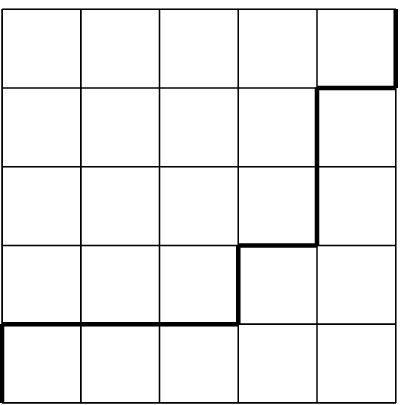}
  \caption{The path $\aPath$.}
  \label{fig:aPath}
  \ifams
\end{figure}
\else
\end{wrapfigure}
\fi
In the following, $P(n,m)$ shall denote the set of words
$w \in \set{x,y}^{\ast}$ such that $\lx{w} = n$ and $\ly{w} = m$.  We
identify a word $w \in P(n,m)$ with a discrete path from $(0,0)$ to
$(n,m)$ which uses only East and North steps of length $1$. For
example, the word $\aPath \in P(5,5)$ is identified with the path in
Figure~\ref{fig:aPath}.

Let $L_{0},L_{1}$ be complete lattices. A \map
$f : L_{0} \rto L_{1} $ is \emph{\jc} if $f(\bigvee X) = \bigvee f(X)$,
for each subset $X$ of $L_{0}$.  We use $\Qj(L_{0},L_{1})$ to denote
the set of \jc \map{s} from $L_{0}$ to $L_{1}$. If
$L_{0} = L_{1} = L$, then we write $\QjL$ for $\Qj(L,L)$.

The set $\Qj(L_{0},L_{1})$ can be ordered pointwise (i.e. $f \leq g$
if and only if $f(x) \leq g(x)$, for each $x \in L_{0}$); with this
ordering it is a complete lattice. Function composition distributes
over (possibly infinite) joins:
\begin{align}
  \label{eq:distr}
  (\bigvee_{j \in J} g_{j}) \circ (\bigvee_{i \in I} f_{i}) & =
  \bigvee_{j \in J, i \in I} (g_{j} \circ f_{i})\,,
\end{align}
whenever $L_{0},L_{1},L_{2}$ are complete lattices,
$\set{f_{i} \mid i \in I} \subseteq \Qj(L_{0},L_{1})$ and
$\set{g_{j} \mid j \in J} \subseteq \Qj(L_{1},L_{2})$.
A \emph{quantale} (see \cite{rosenthal1990}) is a complete lattice
endowed with a semigroup operation $\circ$ satisfying the distributive
law~\eqref{eq:distr}. Thus, $\QjL$ is a quantale, for each complete
lattice $\QjL$.

\bigskip

For $k \geq 0$, we shall use $\I_{k}$ to denote the chain
$\set{0,1,\ldots ,k}$.
Notice that $f : \In \rto \I_{m}$ is \jc if and only if it is monotone
(or order-preserving) and $f(0) = 0$.  For each $n,m \geq 0$, there is
a well-known bijective correspondence between paths in $P(n,m)$ and
\jc \map{s} in $\Qj(\I_{n},\I_{m})$; next, we recall this bijection.
If $w \in P(n,m)$, then the
occurrences of $y$ in $w$ split $w$ into $m+1$ (possibly empty) blocks
of contiguous $x$s, that we index by the numbers $0,\ldots ,m$:
\begin{align*}
  w & = \blockX{w}{0}\cdot y\cdot\blockX{w}{1}\cdot y \ldots
  \blockX{w}{m-1}\cdot y\cdot\blockX{w}{m}\,.
\end{align*}
We call the words
$\blockX{w}{0}, \blockX{w}{1}, \ldots ,\blockX{w}{m} \in
\set{x}^{\ast}$ the $x$-blocks of $w$.
Given $i \in \set{1,\ldots ,n}$, the index of the block of the $i$-th
occurrence of the letter $x$ in $w$ is denoted by
$\blockXno{w}{i}$. We have therefore
$\blockXno{w}{i} \in \set{0,\ldots ,m}$. Notice that $\blockXno{w}{i}$
equals the number of $y$s preceding the $i$-th occurrence of $x$ in
$w$ so, in particular, $\blockXno{w}{i}$ can be interpreted as the
height of the $i$-th occurrence of $x$ when $w$ is considered as a
path.
Similar definitions, $\blockY{w}{j}$ and $\blockYno{w}{j}$, for 
$j = 1,\ldots ,m$, are given for the blocks obtained by splitting $w$
by means of the $x$s:
\begin{align*}
  w & = \blockY{w}{0}\cdot x\cdot\blockY{w}{1}\cdot x \ldots
  \blockY{w}{n-1}\cdot x\cdot\blockY{w}{n}\,.
\end{align*}

The \map $\BLOCKXNO{w}$, sending $i \in \set{1,\ldots ,n}$ to
$\blockXno{w}{i}$, is monotone from the chain $\set{1,\ldots ,n}$ to
the chain $\set{0,1,\ldots ,m}$.  There is an obvious bijective
correspondence from the set of monotone \map{s} from
$\set{1,\ldots ,n}$ to $\I_m = \set{0,1,\ldots ,m}$ to the set
$\Qj(\I_{n},\I_{m})$ obtained by extending a monotone $f$ by setting
$f(0) \eqdef 0$.
We shall tacitly assume this
bijection and, accordingly, we set $\blockXno{w}{0} \eqdef 0$.  Next,
by setting $\blockXno{w}{n+1} \eqdef m$,
we notice that
\begin{align*}
  \len{\blockY{w}{i}} & = \blockXno{w}{i+1} - \blockXno{w}{i}\,,
\end{align*}
for $i = 0,\ldots ,n$, so $w$ is uniquely determined by the \map
$\BLOCKXNO{w}$.  Therefore, the mapping 
sending $w \in P(n,m)$
to $\BLOCKXNO{w}$ is a bijection from $P(n,m)$ to the set
$\Qj(\I_{n},\I_{m})$.
The dominance ordering on $P(n,m)$ arises from the pointwise ordering
on $\Qj(\I_{n},\I_{m})$ via the bijection.

\bigskip

\begin{figure}
  \centering
  \includegraphics{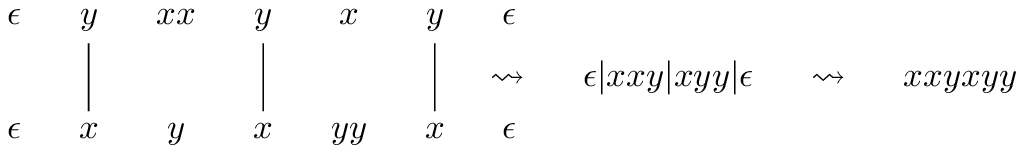}
  \caption{Construction of the product $yxxyxy \otimes xyxyyx$.}
  \label{fig:product}
\end{figure}

For $w \in P(n,m)$ and $u \in P(m,k)$, the product $w \otimes u$ is
defined by concatenating the $x$-blocks of $w$ and the $y$-blocks of
$u$:
\begin{definition}
  For $w \in P(n,m)$ and $u \in P(m,k)$, we let
  \begin{align*}
    w \otimes u & \eqdef \blockX{w}{0} \cdot \blockY{u}{0}\cdot
    \blockX{w}{1}\cdot \blockY{u}{1} \ldots \blockX{w}{m}\cdot
    \blockY{u}{m}\,.
  \end{align*}
\end{definition}
\begin{example}
  Let $w = yxxyxy$ and $u = xyxyyx$, so the $x$-blocks of $w$ are
  $\epsilon,xx,x,\epsilon$ and the $y$-blocks of $u$ are
  $\epsilon,y,yy,\epsilon$; we have $w \otimes u = xxyxyy$.  We can
  trace the original blocks by inserting vertical bars in
  $w \otimes u$ so to separate $\blockX{w}{i}\blockY{u}{i}$ from
  $\blockX{w}{i+1}\blockY{w}{i+1}$, $i = 0,\ldots ,m-1$. That is, we
  can write $w \otimes_{tr} u = |xxy|xyy|$, so $w \otimes u$ is
  obtained from $w \otimes_{tr} u$ by deleting vertical bars. Notice
  that also $w$ and $u$ can be recovered from $w \otimes_{tr} u$, for
  example $w$ is obtained from $w \otimes_{tr} u$ by deleting the
  letter $y$ and then renaming the vertical bars to the letter $y$.
  Figure~\ref{fig:product} suggests that $\otimes$ is a form of
  synchronisation product, obtained by shuffling the $x$-blocks of $w$
  with the $y$-blocks of $u$ so to give ``priority'' to all the $x$s
  (that is, the $x$s precede the $y$s in each block).  It can be
  argued that there are other similar products, for example, the one
  where the $y$s precede the $x$s in each block, so
  $w \oplus u = yxxyyx$.  It is easy to see that
  $w \oplus u = (u^{\star} \otimes w^{\star})^{\star}$, where
  $w^{\star}$ is the image of $w$ along the morphism that exchanges
  the letters $x$ and $y$.
\end{example} 
\begin{proposition}
  The product $\otimes$ corresponds, under the bijection, to function
  composition.  That is, we have
  \begin{align*}
    \BLOCKXNO{w \otimes u} & = \BLOCKXNO{u} \circ \BLOCKXNO{w}\,.
  \end{align*}
\end{proposition}
\begin{myproof}
  In order to count the number of $y$s preceding the $i$-th occurrence
  of $x$ in $w \otimes u$, it is enough to identify the block number
  $j$ of this occurrence in $w$, and then count how many $y$s precede
  the $j$-th occurrence of $x$ in $u$.  That is, we have
  $\blockXno{w \otimes u}{i} = \blockXno{u}{j}$ with $j = \blockXno{w}{i}$.
\end{myproof}
\begin{remark}
Let us exemplify how the algebraic structure of $\Qj(\I_{n},\I_{m})$
yields combinatorial identities.  The product is a function
$\otimes : P(n,m) \times P(m,k) \rto P(n,k)$, so we study how many
preimages a word $w \in P(n,k)$ might have.
By reverting the operational description of the product previously
given, this amounts to inserting $m$ vertical bars marking the
beginning-end of blocks (so to guess a word of the form
$u_{0} \otimes_{tr} u_{1}$) under a constraint that we describe
next. 
Each position can be barred more than once, so adding $j$ bars can be
done in $\binom{n+k + j}{j}$ ways.  The only constraint we need to
satisfy is the following.  Recall that a position
$\ell \in \set{0,\ldots ,n+k}$ is a \emph{\NETurn} (or a
\emph{descent}), see \cite{Krattenthaler97}, if $\ell > 0$,
$w_{\ell-1} = y$ and $w_{\ell} = x$.  If a position is a \NETurn, then
such a position is necessarily barred.  Let us illustrate this with
the word $xxyxyy$ which has just one descent, which is necessarily
barred: $xxy|xyy$.  Assuming $m = 3$, we need to add two more vertical
bars.  For example, for $x|xy|xy|y$ we obtain the following
decomposition:
\begin{align*}
  x|xy|xy|y & \leadsto (x|x|x|,|y|y|y)  \leadsto (xyxyxy,xyxyxy)\,.
\end{align*}
Therefore, if $w$ has $i$ descents, then these positions are barred,
while the other $m -i$ barred positions can be chosen arbitrarily, and
there are $\binom{n+k + m-i}{m- i}$ ways to do this.
Recall that there are $\binom{n}{i}\binom{k}{i}$ words $w \in P(n,k)$
with $i$ descents, since such a $w$ is determined by the subsets of
$\set{1,\ldots ,n}$ and $\set{1,\ldots ,k}$ of cardinality $i$,
determining the descents.  Summing up w.r.t. the number of descents,
we obtain the following formulas:
\begin{align*}
  \binom{n+m}{n}\binom{m+k}{k} & = \sum^{m}_{i = 0} \binom{n+m+k-i}{m-
    i} \binom{n}{i}\binom{k}{i}\,,
  &
  \binom{2n}{n}^{\!2} & = \sum^{n}_{i = 0} \binom{3n-i}{n-
    i}\binom{n}{i}^{\!2}\,.
\end{align*}
Similar kind of combinatorial transformations and identities appear in
\cite{GS78,worpit2011,engbers16}, yet it is not clear to us at the
moment of writing whether these works relate in some way to the
product of paths studied here.
\end{remark}

\begin{remark}
  The previous remark also shows that if $w \in P(n,k)$ has $m \geq 0$
  descents, then there is a canonical factorization
  $w = w_{0} \otimes w_{1}$ with $w_{0} \in P(n,m)$ and
  $w_{1} \in P(m,k)$. It is readily seen that, via the bijection, this
  is the standard epi-mono factorization in the category of
  join-semilattices.  The word $xxyxyy$, barred at its unique descent
  as $xxy|xyy$, is decomposed into $xxyx$ and $yxyy$.
\end{remark}

\begin{remark}
  As in \cite{LaradjiUmar2016}, many semigroup-theoretic properties of
  the monoid $\QjIn$ can be read out of (and computed from) the
  bijection with $P(n,n)$. For example
  \begin{align*}
    \card(\,\set{f \in \QjIn \mid \card(\Image(f)) = k+1}\,) & = \binom{n}{k}^{2}
  \end{align*}
  since, as in the previous remark, a path with $k$ \NETurn{s}
  corresponds to a \jc \map $f$ such that $\card(\Image(f)) = k +1$.
  Similarly
  \begin{align*}
    \card(\,\set{f \in \QjIn \mid \max(\Image(f)) = k}\,) & =
    \binom{n+k-1}{k}
  \end{align*}
  since a \map $f \in \QjIn$ such that $\max(\Image(f)) = k$ (i.e.
  $f(n) = k$) corresponds to a path in $P(n,k)$ whose last step is an
  \EStep, thus to a path in $P(n-1,k)$. A similar argument can be used
  to count \map{s} $f \in \O_{n}$ such that $f(n) = k$,
  cf. \cite[Proposition 3.7]{LaradjiUmar2006}.
\end{remark}

\begin{remark}
  Further properties of the monoid $\QjIn$ can be easily verified, for
  example, this monoid is aperiodic. For the next observation, see
  also \cite[Proposition 2.3]{LaradjiUmar2006} and \cite[Theorem
  3.4]{LaradjiUmar2016}. Recall that $f \in \QjIn$ is \emph{nilpotent}
  if, for some $\ell \geq 0$, $f^{\ell}$ is the bottom of the lattice,
  that is, it is the constant \map with value $0$. It is easily seen
  that $f$ is nilpotent if and only if $f(x) < x$, for each
  $x = 1,\ldots ,n$. Therefore, a path is nilpotent if and only it
  lies below the diagonal, that is, it is a Dyck path. Therefore,
  there are $\frac{1}{n+1}\binom{2n}{n}$ nilpotents in $\QjIn$.
\end{remark}

\section{Idempotent \jc \map{s} as \emmentaler{s}}
\label{sec:emmentalers}

We provide in this section a characterization of idempotent \jc \map{s}
from a complete lattice to itself. 
The characterization originates from the notion of $EA$-duet used to
study some elementary subquotients in the category of lattices, see
\cite[Definition 9.1]{JEMS}.

\begin{definition}
  An \emph{\emmentaler} of a complete lattice $L$ is a collection
  $\E = \set{ [y_{i},x_{i}] \mid i \in I }$ of closed intervals of $L$
  such that
  \begin{itemize}
  \item $[y_{i},x_{i}] \cap [y_{j},x_{j}] = \emptyset$, for $i , j \in
    I$ with $ i \neq j$,
  \item $\set{y_{i} \mid i \in I}$ is a subset of $L$ \caj,
  \item $\set{x_{i} \mid i \in I}$ is a subset of $L$ \cam.
  \end{itemize}
\end{definition}
The main result of this section is the following statement.
\begin{theorem}
  \label{thm:mainEmmentalers}
  For an arbitrary complete lattice $L$, there is a bijection between
  idempotent \jc \map{s} from $L$ to $L$ and \emmentaler{s} of $L$.
\end{theorem}

For an \emmentaler $\E = \set{[y_{i},x_{i}] \mid i \in I}$ of $L$, we
let
\begin{align*}
  \IMG(\E) & \eqdef \set{y_{i } \mid i \in I}\,, & \INCR(\E) &
  \eqdef \set{x_{i } \mid i \in I}\,\\
  \oE(z) & \eqdef \bigvee \set{y \in \IMG(\E) \mid y \leq z}\,, &
  \jE(z) & \eqdef \bigwedge \set{x \in \INCR(\E) \mid z \leq x }\,.
\end{align*}
It is a standard fact that $\jE$ is a closure operator on $L$ (that
is, it is a monotone inflating idempotent \map from $L$ to itself) and
that $\oE$ is an interior operator on $L$ (that is, a monotone,
deflating, and idempotent endomap of $L$).  In the following
statements an \emmentaler $\E = \set{[y_{i},x_{i}] \mid i \in I}$ is
fixed.
\begin{lemma}
  For each $i \in I$, $x_{i} = \jE(y_{i})$ and $\oE(x_{i}) =
  y_{i}$. Therefore $\oE$ restricts to an order isomorphism from
  $\INCR(\E)$ to $\IMG(\E)$ whose inverse is $\jE$ .
\end{lemma}
\begin{myproof}
  Clearly, $\jE(y_{i}) \leq x_{i}$. Let us suppose that
  $y_{i} \leq x_{j}$ yet $x_{i} \not\leq x_{j}$, then
  $y_{i} \leq x_{j} \land x_{i} < x_{i}$ and
  $x_{j} \land x_{i} = x_{\ell}$ for some $\ell \in I$ with
  $\ell \neq i$.  But then
  $x_{\ell} \in [y_{\ell},x_{\ell}] \cap [y_{i},x_{i}]$, a
  contradiction.  
  The equality $\oE(x_{i}) = y_{i}$ is proved similarly.
\end{myproof}

In view of the following lemma we think of $\E$ as a sublattice of $L$
with prescribed holes/fillings, whence the naming ``\emmentaler''.
\begin{lemma}
  If $\E$ is an \emmentaler of $L$, then $\bigcup \E$ is a subset of
  $L$ \cajam. Moreover, the \map sending $z \in [y_{i},x_{i}]$ to
  $y_{i}$ is a complete lattice homomorphism from $\bigcup \E$ to
  $\IMG(\E)$.
\end{lemma}
\begin{myproof}
  Let $\set{z_{k} \mid k \in K}$ with
  $z_{k} \in [y_{k},x_{k}]$ for each $k \in K$. Then, for
  some $j \in I$,
  \begin{align}
    \notag y_{j} &= \bigvee_{k \in K} y_{k} \leq \bigvee_{k \in K}
    z_{k} \leq \bigvee_{k \in K} x_{k} \leq \jE(\bigvee_{k \in
      K} x_{k})\\
    & = \bigvee{\!\!\!}_{\INCR(\E)} \,\set{x_{k} \mid k \in K} =
    \bigvee{\!\!\!}_{\INCR(\E)}\,\set{\jE(y_{k}) \mid k \in K} =
    \jE(\bigvee_{k \in K} y_{k}) = \jE(y_{i}) = x_{j}\,,
    \label{eq:computations}
  \end{align}
  where in the second line we have used the fact that
  $\jE(\bigvee_{k \in K} x_{k})$ is the join in $\INCR(\E)$ of the
  family $\set{x_{k} \mid k \in K}$ and also the fact that $\jE$ is an
  order isomorphism (so it is \jc) from $\IMG(\E)$ to $\INCR(\E)$.
  Therefore,
  $\bigvee_{k \in K} z_{k} \in \bigcup \E$ and, in a similar way,
  $\bigwedge_{k \in K} z_{k} \in \bigcup \E$.

  Next, let $\pi : \bigcup \E \rto \IMG(\E)$ be the \map sending
  $z \in [y_{i},x_{i}]$ to $y_{i} \in \IMG(\E)$. The computations in
  \eqref{eq:computations} show that $\pi$ is \jc. With similar
  computations it is seen that $\bigwedge_{k \in K} z_{k}$ is
  sent to $\oE(\bigwedge_{k \in K} y_{k})$ which
  is 
  the meet of the family $\set{y_{k} \mid k \in K}$ within
  $\IMG(\E)$. Therefore, $\pi$ is \mc as well.
\end{myproof}
We recall next some facts on adjoint pairs of \map{s}, see
e.g. \cite[\S 7]{DP}.  Two monotone \map{s} $f,g : L \rto L$ form an
\emph{adjoint pair} if $f(x) \leq y$ if and only if $x \leq g(y)$, for
each $x,y \in L$. More precisely, $f$ is \emph{left} (or \emph{lower})
\emph{adjoint} to $g$, and $g$ is \emph{right} (or \emph{upper})
\emph{adjoint} to $f$. Each \map determines the other: that is, if $f$
is \ladj to $g$ and $g'$, then $g = g'$; if $g$ is \radj to $f$ and
$f'$, then $f = f'$.  If $L$ is a complete lattice, then a monotone
$f : L \rto L$ is a \ladj (that is, there exists $g$ for which $f$ is
\ladj to $g$) if and only if it is \jc; under the same assumption, a
monotone $g : L \rto L$ is a \radj if and only if it is \mc.

\begin{proposition}
  If $\E$ is an \emmentaler of $L$, then the \map{s} $\fE$ and $\gE$
  defined by
  \begin{align*}
    \fE(z) & \eqdef \oE(\jE(z))\,,
    & 
    \gE(z) & \eqdef \jE(\oE(z))\,,
  \end{align*}
  are idempotent and adjoint to each other. In particular, $\fE$ is
  \jc, so it belongs to $\QjL$.
\end{proposition}
\begin{myproof}
  Clearly, $\fE$ is idempotent:
  \begin{align*}
    \oE(\jE(\oE(\jE(z)))) & = \oE(\jE(z))\,,
  \end{align*}
  since $\jE(z) = x_{i}$ for some $i \in I$ and
  $\jE(\oE(x_{i})) = x_{i}$. In a similar way, $\gE$ is idempotent.
  Let us argue that $\fE$ and $\gE$ are adjoint.  If
  $z_{0} \leq \jE(\oE(z_{1}))$, then
  $\jE(z_{0}) \leq \jE(\jE(\oE(z_{1}))) = \jE(\oE(z_{1}))$ and
  $\oE(\jE(z_{0})) \leq \oE(\jE(\oE(z_{1}))) = \oE(z_{1}) \leq
  z_{1}$. Similarly, if $\oE(\jE(z_{0})) \leq z_{1}$, then
  $z_{0} \leq \jE(\oE(z_{1}))$.
\end{myproof}

\begin{lemma}
  \label{lemma:imagefE}
  $\IMG(\E) = \Image(\fE)$ and $\INCR(\E) = \Image(\gE)$.
\end{lemma}
\begin{myproof}
  Clearly, if $y = \oE(\jE(z))$ for some $z \in L$, then
  $y \in \IMG(\E)$. Conversely, if $y \in \IMG(\E)$, then
  $y = \oE(\jE(y))$, so $y \in \Image(\fE)$.  The other equality is
  proved similarly.
\end{myproof}

For the next proposition, recall that if $f,g$ are adjoint, then
$f \circ g \circ f = f$ and $g \circ f \circ g = g$.
\begin{proposition}
  \label{prop:Ef}
  Let $f \in \QjL$ be idempotent and let $g$ be its \radj. Then
  \begin{enumerate}
  \item $y \leq g(y)$, for each $y \in \Image(f)$,
  \item the collection of intervals
    $\Ef \eqdef\set{[y,g(y)] \mid y \in \Image(f)}$ is an \emmentaler
    of $L$,
  \item   $\IMG(\Ef) = \Image(f)$ and $\INCR(\Ef) = \Image(g)$.
  \end{enumerate}
\end{proposition}
\begin{myproof}
  If $y \in \Image(f)$, then $y = f(y)$ and therefore the relation
  $y \leq g(y)$ follows from $f(y) \leq y$.  The subset $\Image(f)$ is
  \caj since $f$ is \jc. Similarly, $\Image(g)$ is \cam, since $g$ is
  \mc. Let us show that $\set{g(y) \mid y \in \Image(f)} =
  \Image(g)$. To this end, observe that if $x = g(z)$ for some
  $z \in L$, then $x = g(z) = g(f(g(z)))$, so $x = g(y)$ with
  $y = f(g(z))$.

  Finally, let $z \in [y_{1},g(y_{1})] \cap [y_{2},g(y_{2})]$. Then
  $y_{i} = f(y_{i}) \leq f(z) \leq f(g(y_{i}))$.  We already observed
  that $f(g(y_{i})) = y_{i}$, so $y_{i} = f(z)$, for $i = 1,2$. We
  have therefore $y_{1} = y_{2} $ and $g(y_{1}) = g(y_{2})$.
\end{myproof}

\begin{lemma}
  \label{lemma:emmentaler}
  If $f \in \QjL$ is idempotent then, for each $x \in L$,
  \begin{enumerate}
  \item $\oEf(x) \leq f(x)$,
  \item if $f(x) \leq x$, then $f(x) = \oEf(x)$,
  \item if $x \in \INCR(\Ef)$, then $f(x) \leq x$, and so
    $f(x) = \oEf(x)$.
  \end{enumerate}

\end{lemma}
\begin{myproof}
  1. Recall that $\oEf(x) \leq x$ and
  $\oEf(x) \in \IMG(\Ef) = \Image(f)$, so $\oEf(x)$ is a fixed point
  of $f$. Then, using monotonicity, $\oEf(x) = f(\oEf(x)) \leq f(x)$.

  2. From $f(x) \leq x$ and recalling that $\oEf(x)$ is the greatest
  element of $\IMG(\Ef) = \Image(f)$ below $x$,
  it immediately follows that $f(x) \leq \oEf(x)$.
    
  3.  Recall that $\INCR(\Ef) = \Image(g)$, where $g$ is \radj to
  $f$. Let $z$ be such that $x = g(z)$, so we aim at proving that
  $f(g(z)) \leq g(z)$. This is follows from 
  $f(f(g(z))) = f(g(z)) \leq z$ and adjointness.
\end{myproof}
\begin{proposition}
  \label{prop:fEf}
  For each idempotent $f \in \QjL$, we have
  $f = \oEf \circ \jEf = \fEf$. 
\end{proposition}
\begin{myproof}
  Since $\jEf(z) \in \INCR(\Ef)$, then $f(\jEf(z)) = \oEf(\jEf(z))$,
  by the previous Lemma.  Therefore we need to prove that
  $f(\jEf(z)) = f(z)$.  This immediately follows from the relation
  $\jEf = g \circ f$ that we prove next.

  We show that $g(f(z))$ is the least element of $\Image(g)$ above
  $z$.  We have $z \leq g(f(z)) \in \Image(g)$ by adjointness. Suppose
  now that $x \in \Image(g)$ and $z \leq x$. If $y \in L$ is such that
  $x = g(y)$, then $z \leq g(y)$ yields $f(z) \leq y$ and
  $g(f(z)) \leq g(y) = x$.
\end{myproof}

We can now give a proof of the main result of this section,
Theorem~\ref{thm:mainEmmentalers}.
\begin{myproof}[Theorem~\ref{thm:mainEmmentalers}]
  We argue that the mappings $\E \mapsto \fE$ and $f \mapsto \Ef$ are
  inverse to each other.
  
  We have seen (Proposition~\ref{prop:fEf}) that, for an idempotent
  $f \in \QjL$, $\fEf = f$.
  Given an \emmentaler $\E$, we have $\IMG(\E) = \Image(\fE)$ by
  Lemma~\ref{lemma:imagefE}, and $\IMG(\EfE) = \Image(\fE)$, by
  Proposition~\ref{prop:Ef}. 
  Therefore, $\IMG(\E)= \IMG(\EfE)$ and,
  similarly, $\INCR(\E) = \INCR(\EfE)$. Since the two sets $\IMG(\E)$
  and $\INCR(\E)$ completely determine an \emmentaler, we have
  $\E = \EfE$.
\end{myproof}

\section{Idempotent discrete paths}
\label{sec:characterization}

It is easily seen that an \emmentaler of the chain $\In$ can be
described by an alternating sequence of the form
\begin{align*}
 0 = y_{0} \leq x_{0} < y_{1} \leq x_{1} < y_{2} \leq \ldots <
  y_{k} \leq x_{k} = n\,,
\end{align*}
so $\IMG(\E) = \set{0,y_{1},\ldots ,y_{k}}$ and
$\INCR(\E) = \set{x_{1},x_{2},\ldots ,x_{k-1},n}$.  Indeed, $\IMG(\E)$
is \caj if and only if $0 \in \IMG(\E)$, while $\INCR(\E)$ is \cam if
and only if $n \in \INCR(\E)$.

The correspondences between idempotents of $\QjIn$, their paths, and
\emmentaler{s} can be made explicit as follows: for $y \in \IMG(\E)$
such that $y \neq 0$, the path corresponding to $\fE$ touches the
point $(y,y)$ coming from the left of the diagonal; for
$x \in \INCR(\E) \setminus \IMG(\E)$, the path corresponding to $\fE$
touches $(x,x)$ coming from below the diagonal. For
$\E = \set{0< 1 < 2\leq 2 < 3 < 4 }$, with $\IMG(\E) = \set{0,2,3}$
and $\INCR(\E) =\set{1,2,4}$, the path corresponding to $\fE$ is
illustrated in Figure~\ref{fig:idempotentPaths}. On the left of the
figure, points of the form $(x,x)$ with $x \in \INCR(\E)$ are squared,
while points of the form $(y,y)$ with $y \in \IMG(\E)$ are circled.

\newcommand{\mybox}[1][2mm]{\rule{#1}{0mm}\rule{0mm}{#1}}
\begin{figure}
  \centering
  \begin{tikzpicture}
    \diagonal{4};
    \pthNoNETurns{4}{2/2,3/3};
    \node[draw, shape=circle] at (0,0) {}; 
    \node[draw] at (1,1) {\mybox}; 
    \node[draw, shape=circle] at (2,2) {}; 
    \node[draw] at (2,2) {\mybox}; 
    \node[draw, shape=circle] at (3,3) {};
    \node[draw] at (4,4) {\mybox}; 
  \end{tikzpicture}
  \qquad
  \begin{tikzpicture}
    \upDiagonal{4};
    \pthNoNETurns{4}{2/2,3/3};
  \end{tikzpicture}
  \caption{Idempotent path corresponding to $\set{0< 1 < 2\leq 2 <
      3 < 4 }$.}
  \label{fig:idempotentPaths}
\end{figure}
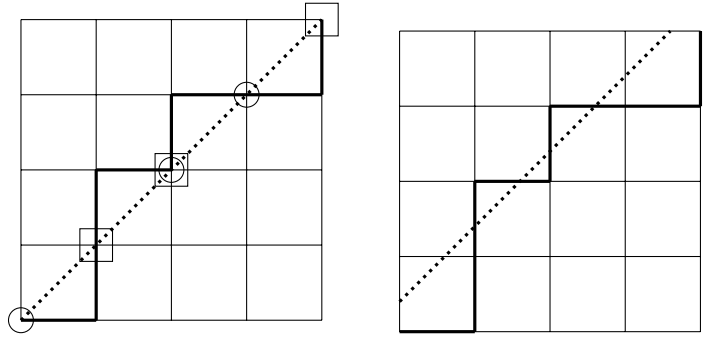

\bigskip

Our next goal is to give a geometric characterization of idempotent
paths using their \NE and \ENTurn{s}.
To this end, observe that we can describe \NETurn{s} of a path
$w \in P(n,m)$ by discrete points in the plane. Namely, if
if $w = w_{0}w_{1} \in P(n,m)$ with $w_{0} = u_{0}y$,
$w_{1} = xu_{1}$, and $\len{w_{0}} = \ell$ (so $w$ has a \NETurn at
position $\ell$), then we can denote this \NETurn with the point
$(\lx{w_{0}},\ly{w_{0}})$.
In a similar way, we can describe \ENTurn{s} by discrete points in the
plane.

Let us call a path an \emph{\uzigzag} if every of its \NETurn{s} is
above the line $y = x + \frac{1}{2}$ while every of its \ENTurn{s} is
below this line. Notice that a path is an \uzigzag if and only every
\NETurn is of the form $(x,y)$ with $x < y$ and every \ENTurn is of
the form $(x,y)$ with $y \leq x$. This property is illustrated on the
right of Figure~\ref{fig:idempotentPaths}.
\begin{theorem}
  \label{thm:geometricCharacterization}
  A path $w \in P(n,n)$ is idempotent if and only if it is an
  \uzigzag.
\end{theorem}

The proof of the theorem is scattered into the next three lemmas.
\begin{lemma}
  An \uzigzag path is idempotent.
\end{lemma}
\begin{myproof}
  Let $w$ be an \uzigzag with
  $\set{(x_{i},y_{i}) \mid i = 1,\ldots ,k}$ the set of its
  \NETurn{s}. For
  $(i,j) \in \set{0,\ldots ,n-1}\times \set{1,\ldots ,n}$, let
  $e_{i,j} \eqdef x^{i}y^{j}x^{n-i}y^{n-j}$ be the path that has a
  unique \NETurn at $(i,j)$.  Notice that
  $w = \bigvee_{i = 1,\ldots ,k} e_{x_{i},y_{i}}$.
  By equation~\eqref{eq:distr},
  \begin{align}
    \label{eq:idempotent}
    w \otimes w & = (\bigvee_{i =1,\ldots ,k} e_{x_{i},y_{i}}) \otimes
    (\bigvee_{j =1,\ldots ,k} e_{x_{j},y_{j}}) = \bigvee_{i,j
      =1,\ldots ,k} e_{x_{i},y_{i}} \otimes e_{x_{j},y_{j}}
    \,.
  \end{align}
  It is now enough to observe that $e_{a,b} \otimes e_{c,d} = e_{a,d}$
  if $c < b$ and, otherwise, $e_{a,b} \otimes e_{c,d} = \bot$, where
  $\bot = x^{n}y^{n}$ is the least element of $P(n,n)$.  Therefore, we
  have: (i)
  $e_{x_{i},y_{i}} \otimes e_{x_{i},y_{i}} = e_{x_{i},y_{i}}$, since
  $x_{i} < y_{i}$, (ii) if $i < j$, then
  $e_{x_{i},y_{i}} \otimes e_{x_{j},y_{j}} = \bot$, since
  $y_{i} \leq x_{j}$, (iii) if $j < i$, then
  $e_{x_{i},y_{i}} \otimes e_{x_{j},y_{j}} = e_{x_{i},y_{j}}$, since
  $x_{j} < y_{j} \leq y_{i}$; in the latter case, we also have
  $e_{x_{i},y_{j}} \leq e_{x_{i},y_{i}}$, since $y_{j} \leq y_{i}$.
  Consequently, the expression on the right of~\eqref{eq:idempotent}
  evaluates to $\bigvee_{i =1,\ldots ,k} e_{x_{i},y_{i}} = w$.
\end{myproof}

Next, let us say that $i \in \In \setminus \set{n}$ is an
\emph{increase} of $f \in \Qj(\In,\I_{m})$ if $f(i) < f(i +1)$.  It is
easy to see that the set of \NETurn{s} of $w$ is the set
$\set{(i,\blockXno{w}{i+1}) \mid i \text{ is an increase of }
  \BLOCKXNO{w}}$.

\begin{lemma}
  \label{lemma:increase}
  Let $f \in \Qj(\I_{n},\I_{n})$ and let $g$ be its right
  adjoint. Then $i \in \I_{n} \setminus \set{n}$ 
  is an increase of $f$ if and only if
  $i \in \Image(g) \setminus \set{n}$.
\end{lemma}
\begin{myproof}
  Suppose $i = g(j)$ for some $j \in \I_{n}$. If $f(i+1) \leq f(i)$,
  then $i + 1 \leq g(f(i)) = g(f(g(j))) = g(j) = i$, a
  contradiction. Therefore $f(i) < f(i+1)$.

  Conversely, if $f(i) < f(i+1)$, then
  $f(i+1) \not\leq f(i)$, $i + 1 \not\leq g(f(i))$, and $g(f(i)) < i +
  1$. Since $i \leq g(f(i))$, then $g(f(i)) = i$, so $i \in \Image(g)$.
\end{myproof}

\begin{lemma}
  The \NETurn{s} of an idempotent path $w \in P(n,n)$ corresponding to
  the \emmentaler
  $\set{0 = y_{0} \leq x_{0} < y_{1},\ldots y_{k} \leq x_{k} = n}$ of
  $\In$ are of the form $(x_{\ell},y_{\ell+1})$, for
  $\ell = 0,\ldots ,k-1$.  Its \ENTurn{s} are of the form
  $(x_{\ell},y_{\ell})$, for $\ell = 0,\ldots ,k$.  Therefore $w$ is
  an \uzigzag.
\end{lemma}
\begin{myproof}
  For the first statement, since
  $\Image(\gE) = \set{x_{0},\ldots ,x_{k-1},n}$ and using
  Lemma~\ref{lemma:increase}, we need to verify that $\fE(x_{\ell}) = y_{\ell}$: this is
  Lemma~\ref{lemma:emmentaler}, point 3.  The last statement is a
  consequence of the fact that \ENTurn{s} are computable from
  \NETurn{s}: if $(x_{i},y_{i})$, $i = 1,\ldots ,k$, are the
  \NETurn{s} of $w$, with $x_{i} < x_{j}$ and $y_{i} < y_{j}$ for
  $i < j$, then \ENTurn{s} of $w$ are of the form $(x_{1},0)$ (if
  $x_{1} > 0$), $(x_{i+1},y_{i})$, $i = 1,\ldots ,k-1$, and
  $(n,y_{k})$ (if $y_{k}< n$).
\end{myproof}

\section{Counting idempotent discrete paths} 

The goal of this section is to exemplify how the characterizations of
idempotent discrete paths given in Section~\ref{sec:characterization}
can be of use.  It is immediate to establish a bijective
correspondence between \emmentaler{s} of the chain $\I_{n}$ and words
$w = w_{0}\ldots w_{n}$ on the alphabet $\set{\mone,0,1}$ that avoid
the pattern $\mone0^{\ast}\mone$ and such that $w_{0} = 1$ and
$w_{n} \in \set{1,\mone}$; this bijection can be exploited for the
sake of counting. We prefer to count idempotents using the
characterization given in
Theorem~\ref{thm:geometricCharacterization}. In the following, we
provide a geometric/combinatorial proof of counting results
\cite{LaradjiUmar2006,howie71} for the number of idempotent elements
in the monoid $\QjIn$ and, also, in the monoid $\O_{n}$ of order
preserving \map{s} from $\set{1,\ldots ,n}$ to itself.
Let us recall that  the Fibonacci
sequence  is defined by 
$f_{0} \eqdef 0$, $f_{1} \eqdef 1$, and
$f_{n+2} \eqdef f_{n +1} + f_{n}$.
Howie \cite{howie71} proved that $\phi_{n} = f_{2n}$ (for $n \geq 1$),
where $\phi_{n}$ is the number of idempotents in the monoid
$\O_{n}$. 
\LU \cite{LaradjiUmar2006} proved that $\gamma_{n} = f_{2n -1}$
($n \geq 1$), where now $\gamma_{n}$ is the number of idempotent
elements of $\O_{n}^{n}$, the submonoid in $\O_{n}$ of \map{s} fixing
$n$. Clearly, $\O_{n}^{n}$ is a monoid isomorphic (and anti-isomorphic
as well) to $\Qj(\I_{n -1})$.
We infer that the number $\psi_{n}$ of idempotents in the monoid
$\QjIn$ equals $f_{2n + 1}$ (for $n \geq 0$).

\begin{remark}
  It is argued in \cite{howie71} that
  $\phi_{n} = \frac{1}{2^{n}\sqrt{5}}\{(3 + \sqrt{5})^{n} - (3 -
  \sqrt{5})^{n}\}$,
  which can easily be verified using the fact that
  $f_{n} = \frac{\theta_{0}^{n} - \theta_{1}^{n}}{\theta_{0} -
    \theta_{1}}$ with $\theta_{0} = \frac{1 + \sqrt{5}}{2}$ and
  $\theta_{1} = \frac{1 - \sqrt{5}}{2}$, see \cite{Gracinda}. In a
  similar way, we derive the following explicit formula: 
  \begin{align*}
    \psi_{n} & = \frac{1}{2^{n+1}\sqrt{5}}\{(3 + \sqrt{5})^{n}(1 +
    \sqrt{5}) - (3 - \sqrt{5})^{n}(1 - \sqrt{5})\} \,.
  \end{align*}
\end{remark}

Let us observe that the monoid $\O_{n}$ can be identified with the
submonoid of $\QjIn$ of \jc maps $f$ such that $1 \leq f(1)$. A path
corresponds to such an $f$ if and only if its first step is a \NStep.
Having observed that $\psi_{0} = \phi_{1} = 1$, the following
proposition suffices to assert that $\phi_{n} = f_{2n}$ and
$\psi_{n} = f_{2n +1}$.
\begin{proposition}
  \label{prop:recurrences}
  The following recursive relations hold:
  \begin{align*}
    \phi_{n+1} & = \psi_{n} + \phi_{n} \,,
    & \psi_{n +1} & = \phi_{n+1} + \psi_{n} \,.
  \end{align*}
\end{proposition}
\begin{myproof}
  Every discrete path from $(0,0)$ to $(n+1,n+1)$ ends with $y$---that
  is, it visits the point $(n+1,n)$---or ends with $x$---that is, it
  visits the point $(n,n+1)$. Consider now an \uzigzag path $\pi$ from
  $(0,0)$ to $(n+1,n+1)$ that visits $(n+1,n)$, see
  Figure~\ref{fig:endinginy}.
  \begin{figure}[h]
    \centering
    \includegraphics[scale=0.7]{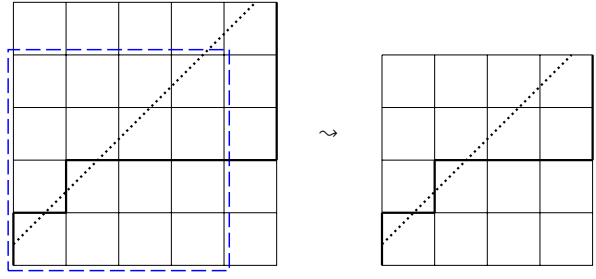}
    \caption{An \uzigzag path to $(5,5)$ ending with $y$.}
    \label{fig:endinginy}
  \end{figure}
  By clipping on the rectangle with left-bottom corner $(0,0)$ and
  right-up corner $(n,n)$, we obtain an \uzigzag path $\pi'$ from
  $(0,0)$ to $(n,n)$. If $\pi$ starts with $y$, then $\pi'$ does as
  well. This proves the right part of the recurrences above, i.e.
  $\phi_{n+1} = \ldots + \phi_{n} $ and
  $\psi_{n +1} = \ldots + \psi_{n}$.

  Consider now an \uzigzag path $\pi$ ending with $x$, see
  Figure~\ref{fig:endinginx}.
  \begin{figure}[h]
    \centering
    \includegraphics[scale=0.7]{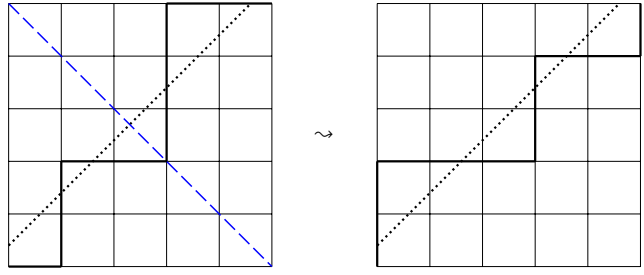}
    \caption{An \uzigzag path to $(5,5)$ ending with $x$.}
    \label{fig:endinginx}
  \end{figure}
  The reflection along the line
  $y = n - x$ sends $(x,y)$ to $(n-y,n-x)$, so it preserves \uzigzag
  paths. Applying this reflection to $\pi$, we obtain an \uzigzag path
  from $(0,0)$ to $(n+1,n+1)$ whose first step is $y$.  This proves
  the $\psi_{n+1} = \phi_{n+1} + \ldots $ part of the recurrences
  above.

  Consider now an \uzigzag path $\pi$ ending with $x$ and beginning
  with $y$, see Figure~\ref{fig:endinginxbegy}.
  \begin{figure}[h]
    \centering
    \includegraphics[scale=0.7]{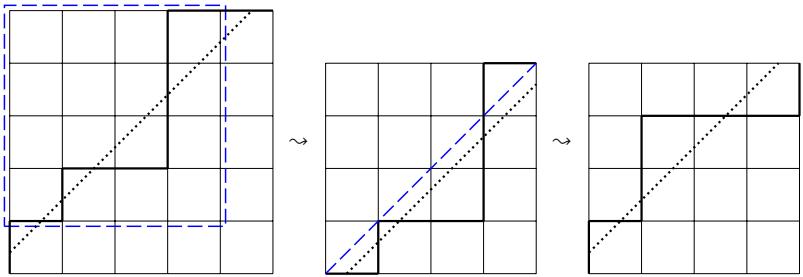}
    \caption{An \uzigzag path to $(5,5)$ ending with $x$ and beginning
      with $y$.}
    \label{fig:endinginxbegy}
  \end{figure}
  By clipping on the
  rectangle with left-bottom corner $(0,1)$ and right-up corner
  $(n,n+1)$ and then by applying the translation $x \mapsto x -1$, we
  obtain a path whose all \NETurn{s} are above the line
  $y = x - \frac{1}{2}$ and whose all \ENTurn{s} are below this
  line. By reflecting along diagonal, we obtain an \uzigzag path from
  $(0,0)$ to $(n,n)$. This proves the
  $\phi_{n+1} = \psi_{n} + \ldots $ part of the recurrences above.
\end{myproof}
The geometric ideas used in the proof of
Proposition~\ref{prop:recurrences} can be exploited further, so to
show that the number of idempotent maps $f \in \Qj(\In)$ such that
$f(n) = k$ equals $f_{2k}$, see the analogous statement in
\cite[Corollary 4.5]{LaradjiUmar2006}.  Indeed, if $f(n) = k$, then
the path corresponding to $f$ visits the points $(n-1,k)$ and $(n,k)$;
therefore, since it is an \uzigzag, also the points $(k-1,k)$ and
$(k,k)$. By clipping on the rectangle from $(0,0)$ to $(k,k)$, we
obtain an \uzigzag path in $P(k,k)$ ending in $x$. As seen in the
proof of Proposition~\ref{prop:recurrences}, these paths bijectively
correspond to \uzigzag paths in $P(k,k)$ beginning with $y$.

\section{Conclusions}

We have presented the monoid structure on the set $P(n,n)$ of discrete
lattice paths (with North and \ESteps) that corresponds, under a
well-known bijection, to the monoid $\QjIn$ of \jc functions from the
chain $\set{0,1,\ldots n}$ to itself.  In particular, we have studied
the idempotents of this monoid, relying on a general characterization
of idempotent \jc functions from a complete lattice to itself.  This
general characterization yields a bijection with a language of words
on a three letter alphabet and a geometric description of idempotent
paths. Using this characterization, we have given a
geometric/combinatorial proof of counting results for idempotents in
monoids of monotone endo\map{s} of a chain
\cite{howie71,LaradjiUmar2006}.

Our initial motivations for studying idempotents in $\Qj(\In)$
originates from the algebra of logic, see e.g. \cite{Jipsen}. Willing
to investigate congruences of $\Qj(\In)$ as a residuated lattice
\cite{GJKO},
it can be shown, using idempotents, that every subalgebra of a
residuated lattice $\QjIn$ is simple. This property does not
generalize to infinite complete chains: if $\I$ is the interval
$[0,1] \subseteq \mathbb{R}$, then $\Qj(\I)$ is simple but has
subalgebras that are not simple \cite{BallDroste1985}.
Despite the results we presented are not related to our original
motivations, we aimed at exemplifying how a combinatorial approach
based on paths might be fruitful when investigating various kinds of
monotone \map{s} and the multiple algebraic structures these maps may
carry.
 
We used the Online Encyclopedia of Integer Sequences to trace related
research. In particular, we discovered Howie's work \cite{howie71} on
the monoid $\O_{n}$ through the OEIS sequences 
\OEIS{A001906} 
and \OEIS{A088305}. 
The sequence $\psi_{n}$ is a shift of 
the sequence \OEIS{A001519}. Related to this sequence is the doubly
parametrized sequence \OEIS{A144224} collecting some counting results
from \cite{LaradjiUmar2006} on idempotents.
Relations with other kind of combinatorial objects counted by the
sequence $\psi_{n}$ still need to be understood.

\ifams
\bigskip
\fi

\paragraph{Acknowledgment.} 
The author is thankful to 
Srecko Brlek, Claudia Muresan, and Andr\'e Joyal for the fruitful
discussions he shared with them
on this topic during winter 2018.
The author is also thankful to the anonymous referees for their
insightful comments and for pointing him to the reference
\cite{LaradjiUmar2006}.

\ifams\else
\newpage
\fi

\newpage
\bibliographystyle{splncs04}
\bibliography{biblio}

\end{document}